\documentclass{amsart}
\usepackage{amsfonts}

\begin{document}

\title[COMPLEX-VALUED BEST LINEAR UNBIASED ESTIMATOR OF ...]
{COMPLEX-VALUED BEST LINEAR UNBIASED ESTIMATOR OF AN UNKNOWN
CONSTANT MEAN OF WHITE NOISE}

\author{T. SUS{\L}O}
\email{tomasz.suslo@gmail.com}
\dedicatory{to Els Van Hecke}

\begin{abstract}
In this paper the complex-valued 
best linear unbiased estimator of an unknown constant mean
of white noise was derived  
the ordinary least-squares estimator of an unknown constant 
mean of random field (arithmetic mean) charged by an imaginary error. 
\end{abstract}

\maketitle
 
\thispagestyle{empty}

\section{Introduction}
\noindent
Let us consider a stationary random process 
$\epsilon=\{\epsilon_j=\epsilon(x_j);~x_j \supset x_i=x_1,\ldots,x_n\}$ 
with zero mean 
$$
E\{\epsilon_j\}=E\{\epsilon_i\}=E\{\epsilon\}=0 
$$ 
and the background trend 
$
\sum_{k}f_{jk} \beta^k
=
f_{jk} \beta^k
=
m(x_j)
$ 
(some known mean function $m(x_j)$ with unknown regression 
parameters $\beta^k$) then 
\begin{equation}
V_i
=
\epsilon_i+m(x_i) 
=
\epsilon_i+f_{ik} \beta^k 
\qquad \mbox{for} \quad x_i=x_1,\ldots,x_n 
\label{1}
\end{equation}
and
\begin{equation}
V_j
=
\epsilon_j+m(x_j) 
=
\epsilon_j+ f_{jk}\beta^k 
\label{2}
\end{equation}
where $f_{jk}$ 
is a given vector and $f_{ik}$ 
is a given matrix.

\noindent
The unbiasedness constraint 
$$
E\{V_j\}=E\{\hat{V}_j\}
$$
on the estimation statistics
$$
\hat{V}_j 
=
\omega^i_j V_i
=
\omega^i_j \epsilon_i + \omega^i_j f_{ik} \beta^k
=
\hat{\epsilon}_j + \omega^i_j f_{ik} \beta^k
$$   
produces the system of $N(k)$ equations in $n$ unknowns $\omega^i_j$  
\begin{equation}
f_{jk}=\omega^i_j f_{ik} 
\label{us}
\end{equation}
and gives
\begin{equation}
\hat{V}_j 
=
\hat{\epsilon}_j + f_{jk} \beta^k \ .
\label{3}
\end{equation}

\noindent
For white noise  
\begin{eqnarray}
E\{[V_j-\omega^i_j V_i]^2\}
&
=
&
E\{[V_j-\hat{V}_j]^2\} \nonumber
\\
&
=
&
E\{[\epsilon_j-\hat{\epsilon}_j]^2\} \nonumber 
\\
&
=
&
E\{\epsilon_j^2\}
+
E\{\hat{\epsilon}_j^2\} \nonumber 
\\
&
=
&
E\{[V_j-f_{jk}\beta^k]^2\}
+
E\{[\hat{V_j}-f_{jk} \beta^k]^2\} \nonumber
\\
&
=
&
\sigma^2+\sigma^2 \omega^i_j \rho_{il} \omega^l_j       
\label{oo}
\end{eqnarray}
where $\rho_{il};~i,l=1,\ldots,n$ 
is the indentity matrix.  

Our aim is to constraint  from~(\ref{oo})  
the estimation statistics of the field $V_j$ 
$$
\hat{V}_j
=
\omega^i_j V_i
=
\omega^i_j \epsilon_i + \omega^i_j f_{ik} \beta^k
=
\hat{\epsilon}_j + f_{jk} \beta^k
$$ 
for the estimation statistics of  
mean $f_{jk}\beta^k$ of the field 
at some $x_j$    
given by the constraint on $x_j$
\begin{equation}
E\{[\omega^i_j V_i-f_{jk} \beta^k]^2\}
=
\sigma^2 \omega^i_j \rho_{il} \omega^l_j 
=
0       
\label{c}
\end{equation}
then at some $x_j$ holds
$$
E\{[V_j-\omega^i_j V_i]^2\}
=
\sigma^2
=
E\{[V_j-f_{jk}\beta^k]^2\} \ .               
$$

\section{The best linear unbiased estimation statistics}
\noindent
The minimization constraint on $E\{[\omega^i_j V_i-V_j]^2\}$ 
$$
\frac{\partial 
E\{[\omega^i_j V_i-V_j]^2\}
}{\partial \omega^i_j} 
=
2\sigma^2\rho_{il}\omega^l_j 
+ 2\sigma^2 f_{ik}\mu^k_j = 0 \ ,
$$
where~(\ref{oo})
$$
E\{[\omega^i_j V_i-V_j]^2\}
=\sigma^2 + \sigma^2 \omega^i_j \rho_{il} \omega^l_j  
+ 2\sigma^2\underbrace{(\omega^i_j f_{ik} - f_{jk})}_0 \mu^k_j \ ,
$$
let us add $n$ equations in $N(k)$ unknowns $\mu_j^k$ to 
the unbiased system~(\ref{us}) of $N(k)$ equations in $n$ unknowns $\omega^i_j$  
\begin{equation}
\rho_{is}\omega^s_j  = - f_{ik} \mu^k_j \ , \quad s=1,\ldots,n \ ,
\label{Lw}
\end{equation}
then
\begin{equation}
\delta^l_s\omega^s_j  
=
\omega^l_j
=
- \rho^{li} f_{ik} \mu^k_j \ ,  
\label{w}
\end{equation}
where
$$
\rho^{li}\rho_{is}= \delta^l_s \ ,
$$
substituted into the unbiased system~(\ref{us}) 
$$
f_{zj}=f_{zl} \omega^l_j \ ,\quad z=1,\ldots, N(k) \ ,
$$
gives
\begin{equation}
\mu^k_j   
=
-(f_{zl} \rho^{li} f_{ik})^{-1} f_{zj}   
\label{u}
\end{equation}
and~(from(\ref{w})) the kriging weights
$$
\omega^l_j  
= 
\rho^{li} f_{ik} (f_{zl}\rho^{li} f_{ik})^{-1} f_{zj} \ .
$$
Now, we can derive the kriging estimator  
\begin{equation}
\hat{v}_j
=
v_l \omega^l_j 
=
\omega^l_j v_l 
=  
f_{jz}\hat{\beta}^z 
\label{ke}
\end{equation}
where the ordinary least-squares estimator
\begin{equation}
\hat{\beta}^z=(f_{ki}\rho^{il}f_{lz})^{-1}f_{ki}\rho^{il}v_l 
\label{b}
\end{equation}
is the best linear unbiased estimator of 
unknown regression parameters $\beta^z$ 
based on $v_l$ as an observation of a stochastic
process $V=\{V_j=V(x_j);~x_j \supset x_l=x_1,\ldots,x_n\}$ 
and minimized variance of the best linear unbiased ordinary
(estimation) statistics $\hat{V}_j=\omega^i_j V_i$ 
of the field 
$V_j$~(from~(\ref{oo})~and~(\ref{Lw})~and~(\ref{us})~and~(\ref{u})) 
\begin{eqnarray}
E\{\hat{\epsilon}_j^2\}
=
E\{[\omega^i_j V_i-f_{jk} \beta^k]^2\}
&
=
&
\sigma^2 \omega^i_j \rho_{il} \omega^l_j \nonumber \\
& 
=
&
-\sigma^2 \omega^i_j f_{ik} \mu^k_j \nonumber \\   
&
=
&
-\sigma^2 f_{jk}\mu^k_j \nonumber \\ 
&
=
&
\sigma^2 f_{jk} (f_{zl} \rho^{li} f_{ik})^{-1} f_{zj}  
\label{mvoes} \ .
\end{eqnarray}

\section{Complex-valued estimator of an unknown constant mean of white noise}

\noindent
Since for constant mean function   
$$
V_i=\epsilon_i+f_{ki} \beta^k=\epsilon_i+\beta^{\it 1} 
\qquad \mbox{for} \quad x_i=x_1,\ldots,x_n   
$$
and
$$
V_j=\epsilon_j+f_{kj} \beta^k=\epsilon_j+\beta^{\it 1}  
$$
the minimized variance~(\ref{mvoes}) of the best linear unbiased ordinary
(estimation) statistics
$\hat{V}_j=\omega^i_j V_i$ of the field $V_j$ 
can not be compared to zero value 
$$
E\{\hat{\epsilon}_j^2\}
=
E\{[\hat{V}_j-\beta^{\it 1}]^2\}
=
E\{[\omega^i_j V_i-\beta^{\it 1}]^2\}
=
n^{-1}\sigma^2 
$$
let us consider linear mean function
with the slope $\beta^{\it 2}$ and the offset $\beta^{\it 1}$  
$$
V_i=\epsilon_i+f_{ik} \beta^k=\epsilon_i+\beta^{\it 1}+x_i\beta^{\it 2} 
\qquad \mbox{for} \quad x_i=x_1,\ldots,x_n   
$$
and
$$
V_j=\epsilon_j+f_{jk} \beta^k=\epsilon_j+\beta^{\it 1}+x_j\beta^{\it 2}    
$$
to constraint the estimation statistics 
of the field $V_j$  
$$
\hat{V}_j=\omega^i_j V_i
=
\hat{\epsilon}_j + \beta^{\it 1}+x_j \beta^{\it 2} 
$$ 
for the best linear unbiased statistics of mean  
$f_{jk}\beta^k=\beta^{\it 1}+x_j \beta^{\it 2}$ of the field at some $x_j$ 
we have to solve~(from~(\ref{c})~and~(\ref{mvoes}))
$$
f_{jk} 
(f_{zl}\rho^{li} f_{ik})^{-1}
f_{zj}= 
\underbrace{
\left[
\begin{array}{cc}
1 & x_j \\
\end{array}    
\right]}_{1 \times 2}
{\underbrace{
\left[
\begin{array}{cc}
n & n\overline{x_i} \\
  &            \\    
n\overline{x_i} & n\overline{x_i^2} \\
\end{array}    
\right]}_{2 \times 2}}^{-1}
\underbrace{
\left[
\begin{array}{c}
1 \\
x_j \\
\end{array}    
\right]}_{2 \times 1}
=
\frac
{x_j^2-2m_n x_j+m_{sn}}
{n\sigma^2_n}
=
0 \ ,  
$$
where 
$$
m_n=\overline{x_i}=\frac{1}{n} \sum_i x_i,
\quad 
m_{sn}=\overline{x_i^2}=\frac{1}{n} \sum_i x_i^2,
\quad
\sigma_n=\sqrt{\overline{x_i^2}-{\overline{x_i}}^2} 
$$
and we get
\begin{equation}
x_j = m_n \pm I\sigma_n \ , 
\label{j}
\end{equation}
where $I=\sqrt{-1}$. 

The best linear unbiased estimator 
$$
\hat{v}_j
=
\hat{\beta}^{\it 1} + x_j \hat{\beta}^{\it 2} 
$$
with the ordinary least-squares estimator of the offset $\beta^{\it 1}$
$$
\hat{\beta}^{\it 1}
=
\frac{m_{sn}\overline{v_i}-m_n\overline{x_iv_i}}{\sigma_n^2}
$$
the ordinary least-squares estimator of the slope $\beta^{\it 2}$
$$
\hat{\beta}^{\it 2}
=
\frac{\overline{x_iv_i}- m_n \overline{v_i}}{\sigma_n^2} \ , 
$$ 
based on $v_i=v_1,\ldots,v_n$, where 
$$
\overline{x_i v_i}=\frac{1}{n}\sum_i x_i v_i \ , 
\quad
\overline{v_i}=
\frac{1}{n} \sum_i v_i \ ,
$$
constrained at~(\ref{j})  
$$
x_j = m_n \pm I\sigma_n 
$$ 
is the ordinary least-squares estimator of an unknown constant mean 
of the field (arithmetic mean)
$$
\overline{v_i} 
=
\frac{1}{n} \sum_i v_i 
$$
charged by an imaginary error (Van Hecke estimator)
$$
\hat{v}_j
=
\hat{\beta}^{\it 1} + m_n \hat{\beta}^{\it 2} \pm I \sigma_n \hat{\beta}^{\it 
2}
=
\overline{v_i} \pm I \sigma_n \hat{\beta}^{\it 2}=\hat{m} \ . 
$$

\vspace*{12pt}
\noindent
{\bf Example.} The best linear unbiased estimator of an unknown constant mean  

$$
\overline{v_i}
=
\Re(\hat{m})
=
\frac{1}{n} \sum_i v_i 
= 
3.29
$$
of an uncorrelated signal (white noise) $v_i$ at $x_i$   
$$
\begin{array}{|c|c|c|c|c|c|c|c|c|c|c|c|}
\hline
x_i & 1 & 2 & 3 & 4 & 5 & 6 & 7 & 8 & 9 & 10 & 11 \\     
\hline
v_i & 4.12 & 1.38 & 5.71 & 1.25 & 2.24 & 0.81 & 1.67 & 7.42 & 7.91 & 1.63 & 2.05 \\
\hline
\end{array}
$$
has the real-valued standard error  
$$
\pm \sqrt{\frac{\overline{v_i^2}-\overline{v_i}^2}{n}} = \pm 0.74 
$$
and the imaginary error 
$$
\Im(\hat{m})
=
\pm\sigma_n
\hat{\beta}^{\it 2} 
=
\pm\frac{\overline{x_iv_i}-\overline{x_i}~\overline{v_i}}
{\sqrt{\overline{x_i^2}-{\overline{x_i}}^2}}
= 
\pm 0.26 \ .
$$

\end{document}